\documentclass[12pt]{amsart}
\usepackage{amsfonts}
\usepackage{amsmath}
\usepackage{amsxtra}
\usepackage{amssymb,latexsym}
\usepackage[mathcal]{eucal}
\usepackage{amscd}
\newtheorem{theo}{{\bfseries Theorem}}[section]
\newtheorem{prop}[theo]{{\bfseries Proposition}}
\newtheorem{lem}[theo]{{\bfseries Lemma}}
\newtheorem{cor}[theo]{{\bfseries Corollary}}
\newtheorem{df}[theo]{{\bfseries Definition}}

\def \R {\mathbb R}

\def \A {\mathcal A}
\def \B {\mathcal B}

\def \L {\mathcal L}
\def \J {\mathcal J}

\def \G {\mathcal G}

\def \ep {\epsilon}
\def \om {\omega}
\def \d {\delta}

\def \i {\iota}

\def \t {\tau}

\usepackage{amssymb,latexsym}
\usepackage[mathcal]{eucal}
\usepackage{amscd}
\usepackage{amsmath}
\usepackage{graphicx}
\usepackage{amsfonts}
\usepackage{amscd}

\numberwithin{equation}{section}

\begin{document}

\title[Induced Fuzzy Topological Spaces]{\bfseries  Induced Fuzzy Topological Spaces: \\
A Characterization}
\vspace{1cm}
\author{Ethan Akin}
\address{Mathematics Department \\
    The City College \\ 137 Street and Convent Avenue \\
       New York City, NY 10031, USA     }
\email{ethanakin@earthlink.net}

\date{September, 2017}

\begin{abstract} We introduce a simple property, \emph{affine invariance}, which characterizes within the class of fuzzy topological spaces
those which are induced from an underlying topology on the space.  We illustrate it by considering the simple notions of compactness for such spaces.
\end{abstract}

\keywords{fuzzy topological space, induced fuzzy topological space, affinely invariant fuzzy topologies, fuzzy Tychonoff theorem.}

\vspace{.5cm} \maketitle

\emph{In Memory of Jack Feldman}

%\tableofcontents

\section*{Introduction}

The concept of a fuzzy set was introduced in 1965 by Zadeh \cite{Z65} and independently, as is pointed out by \v{S}ostak \cite{S96}, by
Salij \cite{Sa65}.
A fuzzy subset of a set $X$ is a function $f : X \to I$ where $I = [0,1]$ is
the unit interval in $\R$.  For a point $x \in X$ the value $f(x)$ is understood as the degree of membership
of $x$ in the fuzzy set. An ordinary subset $A$ of $X$ corresponds to its characteristic function $\chi_A$
which is $1$ on $A$ and $0$ on its complement.

The set $I^X$ of functions from $X$ to $I$, i.e. the set of fuzzy subsets, is a complete lattice.
For any family $\{ f_i \in I^X : i \in K \}$ the
pointwise supremum $\bigvee_{i \in K} \ f_i$ and pointwise infimum $\bigwedge_{i \in K} \ f_i$ lie in $I^X$.
The association $A \mapsto \chi_A$ includes the lattice of subsets of $X$ into the lattice $I^X$.

The concept of a fuzzy topology was defined by  Chang \cite{Ch68} as a collection $\d$ of fuzzy sets which satisfies
axioms strictly analogous to the axioms for a topology. The constants $0, 1$ are in $\d$ and $\d$ is closed under arbitrary
supremum and finite infimum.

 Fuzzy topological spaces quickly became the object of study
by a number of authors.  Of particular interest is a series of beautiful papers by Lowen, \cite{L76}, \cite{L77}, \cite{L78}, and \cite{L81}.
The subject has grown immensely.  The current state is  summarized in detail in Palaniappan's book \cite{P02}.
The magisterial survey by \v{S}ostak \cite{S96} covers not only fuzzy topological spaces, but also the broader theory, initiated by Goguen \cite{G67}
(and Salij \cite{Sa65}) where the fuzzy objects take values in a general complete, distributive lattice instead of $I$.

From the beginning there has been some dispute about how broad the definition of a fuzzy topology should be.  In \cite{L76} Lowen introduced among the
axioms a simple additional condition, that the constants $k$ be elements of $\d$ for all $k \in I$.
He provided a number of justifications for his view that the concept of
fuzzy topological space should be restricted to include this strengthened axiom.  These arguments were vigorously extended in \cite{LW88}.
Thus, Lowen and Wuyts use the label quasi-fuzzy topology for the original Chang axioms, reserving fuzzy topology for those which include constants.
On the other hand, \v{S}ostak retains the original Chang definition and uses the term laminated fuzzy topology for one which includes constants.
He argues,\cite{S96} page 671:

\begin{quote} R. Lowen and P. Wuyts insist on restricting the subject of fuzzy topology
entirely to laminated spaces. However, at present, the attitude of
most workers in fuzzy topology towards the property of laminatedness can be compared with
the attitude of general topologists towards the axiom, say, of
Hausdorffness or complete regularity: One may certainly accept it when it is really needed,
but it is not reasonable to include axiom $T_2$ or $T_{3.5}$ as part of the general definition
of a topological space. Thus the standard framework for most specialists working in fuzzy topology is the category
of Chang Fuzzy Topological Spaces, but when it is indeed necessary, they impose the additional
condition of laminatedness on the considered spaces. \end{quote}

I tend to agree with Lowen and Wuyts. In fact, I believe that the subject should be further restricted by the introduction of what I call
affine invariance. However, to avoid confusion, I will follow \v{S}ostak's terminology.

In the early paper \cite{L76} Lowen introduced functors between topological spaces and fuzzy topological spaces.
He pointed out that for a topological space $(X,\t)$ the set of lower
semicontinuous functions in $I^X$ is a laminated fuzzy topology $\om(\t)$. On the other hand, given a
fuzzy topology $\d$ we can associate to it the coarsest topology $\i(\d)$ on $X$ with respect to which the
the functions in $\d$ are lower semicontinuous, that is, the topology with subbase $\{ f^{-1}(c,\infty) :
c \in I $ and $ f \in \d \}$. Lowen also introduced the appropriate notion of compactness for fuzzy
subsets of a space with a laminated fuzzy topology. Not every fuzzy topology, not even every laminated fuzzy topology,
is of the form $\om(\t)$ for some topology $\t$. The fuzzy topologies of this form are said to be topologically induced or simply
induced fuzzy topologies.

We will say that a fuzzy topology $\d$ is affine invariant if for all $m, k \in \R$ with $m > 0$,  and all $f \in \d$ such that
$m \cdot f  +  k \in I^X$, i.e. $m \cdot f(x)  +  k \in I$ for all $x \in X$, the affine adjustment $m \cdot f  +  k $ is an element of $ \d$.  That is,
we assume that $\d$ is closed under appropriate positive affine transformations.
Clearly, the set of lower semicontinuous functions on a topological space is affine invariant. We show, conversely, that
the induced fuzzy topologies are exactly the affine invariant fuzzy topologies.

Finally, we illustrate the utility of the restriction by considering the meaning of compactness in the context of affine invariant fuzzy topologies.

The work here was inspired by a reading of the lovely survey, \cite{C08}. Although my intent is to
contract the range of fuzzy topology, I nonetheless appreciate Carlson's more expansive exposition.
After reading Carlson's paper several summers ago, I had some lunchtime conversations with
Jack Feldman from which this work is the belated fruit.
I dedicate it to his memory.
\vspace{1cm}

\section{Affine Invariance}

Let $(X, \t)$ be a topological space.  A function $f : X \to \R$ is \emph{lower semicontinuous}
(hereafter  \emph{lsc}) if for every $c \in \R$ the set $f^{-1}(c,\infty) \in \t$, i.e. it  is open in $X$.
We will be using this sort of set so often that we introduce the notation:
\begin{equation}\label{01}
f_{(c} \ =_{def} \ f^{-1}(c,\infty) \qquad \mbox{and} \qquad f_{[c} \ =_{def} \ f^{-1}[c,\infty),
\end{equation}
for any $c \in \R$.

If $U$ is a  subset of $X$ then the characteristic function $\chi_U$ of $U$ is lsc if and only if $U \in \t$. Recall
that $\chi_U(x) = 1$ if $x \in U$ and $= 0$ otherwise.

For any  $L \subset \R$ we will denote by $L^X$ the set of all functions from $X$ to $L$.

\begin{df}\label{defa01} For a set $X$ let $\d \subset I^X$. We call $\d$ a \emph{fuzzy topology }
if it satisfies the conditions:
\begin{itemize}
\item[(i)]  The constant functions $0, 1$ are in $\d$.
\item[(ii)] If $\{ f_i : i \in K \}$ is any indexed family of members of $\d$ then
$\bigvee_{i \in K} \ f_i \ \in \ \d$.
\item[(iii)] If $\{ f_i : i \in K \}$ is any indexed family of members of $\d$  with $K$ finite then
$\bigwedge_{i \in K} \ f_i \ \in \ \d$.
\end{itemize}
A pair $(X,\d)$ where $X$ is a set and $\d$ is a fuzzy topology on $X$ will be called a \emph{fuzzy topological space}.

We call $\d$ a \emph{laminated fuzzy topology } (or we say that $(X, \d)$ is a \emph{laminated fuzzy topological space })
when (ii) and (iii) hold and (i) is strengthened to:
\begin{itemize}
\item[(i')] The constant functions $k $ are in $\d$ for all $k \in I$.
\end{itemize}

We call $\d$ an \emph{affine invariant fuzzy topology} (or we say that the fuzzy topological space
$(X, \d)$ is \emph{affine invariant})
when (ii) and (iii) hold and (i) is replaced by:
\begin{itemize}
\item[(i'')] The constant function $0$ is in $ \d$ and $\d$ is closed under suitable positive affine transformations, i.e.
if $f \in \d$, $m,k \in \R$ with $m > 0$ and $g = m \cdot f + k$ satisfies $g(x) \in I$ for all $x \in X$
then $g \in \d$.
\end{itemize}
\end{df}
\vspace{.5cm}

For an indexed family we will always assume that the index set $K$ is  nonempty.

Clearly, these are successively stronger conditions. In particular, (i'') implies (i') [Let $f = 0$].

The concept of fuzzy topology is due to Chang \cite{Ch68}. Lowen introduced the sharper condition (i') in
\cite{L76}. We adopt the term \emph{lamination} following \cite{S96}.

As with topologies, the intersection of arbitrary families of fuzzy topologies on $X$ is a fuzzy topology on $X$.
Conditions (i') and (i'') are preserved by intersection as well.  Hence, for any  $\A \subset I^X$ we can obtain
the smallest fuzzy topology on $X$ which contains $\A$, or the \emph{fuzzy topology generated by $\A$}, by intersecting
all fuzzy topologies which contain $\A$. Notice that $I^X$ is itself an affine invariant fuzzy topology.
It is called the \emph{discrete fuzzy topology}.

Following \cite{W75} and \cite{L76} we introduce constructions relating topologies and fuzzy topologies. We let $\R_r$ denote $\R$ equipped with
the \emph{lower topology} whose non-empty open sets are $\{ (c,\infty) : c \in \R \}$ and we let $I_r$ denote set $I$ with the
relative topology induced from $\R_r$. Thus,
a function $f : X \to \R_r$ is continuous exactly when it is a lower semicontinuous real-valued (hereafter, lsc) function.

\begin{df}\label{defa02} For $(X, \t)$ any topological space let $\om(\t)$ denote the set of all continuous functions from $X$ to $I_r$, or, equivalently,
the set of lsc functions from $X$ to $I$.
That is, $f \in \om(\t)$ if and only if $f \in I^X$ with $f_{(c} \in \t$ for all $c \in I$.

For $(X,\d)$ any fuzzy topological space let $\i(\d)$ denote the coarsest topology on $X$ so that every $f \in \d$ is continuous
from $X$ to $I_r$.  That is, $\i(\d)$ has subbase
$\{ f_{(c} : \ c \in I \ $ and $ \ f \in \d \ \}$.

If $\d$ is a fuzzy topology on $X$ we will say that $\d$ is \emph{induced from the topology $\t$} on $X$ when $\d = \om(\t)$.
We will say that $\d$ is an \emph{induced fuzzy topology} or a \emph{topological fuzzy topology} when it is induced from
some topology on $X$.
\end{df}
\vspace{.5cm}

\begin{theo}\label{theoa03} If $\t$ is a topology on $X$, then $\om(\t)$ is an affine invariant fuzzy topology and
$ \i(\om(\t)) \ = \ \t$. Thus, if $\d$ is the fuzzy topology induced from $\t$, then $\t = \i(\d)$.
\end{theo}

{\bfseries Proof:}  It is easy to check that the set of lsc functions satisfies (ii)and (iii). In general, a continuous
function $\phi : \R_r \to \R_r$ is a non-decreasing function on $\R$ which is continuous from the left. For example,
if $\phi(t) = m \cdot t + b$ with $m \geq 0$ then $\phi$ is such a function. If $f : X \to I_r$ is continuous and
$\phi : \R_r \to \R_r$  is  continuous with $\phi(f(X)) \subset I$ then $\phi \circ f : X \to I_r$ is continuous.  Thus,
$f \in \om(\t)$ implies that $\phi \circ f \in \om(\t)$. Hence,
 $\om(\t)$ is an affine invariant fuzzy topology.

By definition if $f$ is lsc on $(X,\t)$ then each $f_{(c}$ is in $\t$ and so
the inclusion $\i(\om(\t)) \subset \t$ is clear.  On the other hand, if $U \in \t$ then $\chi_U \in \om(\t)$ and
so $U \in \i(\om(\t))$.

Thus, if $\d = \om(\t)$ then $\i(\d) = \i(\om(\t)) \ = \ \t$.

$\Box$ \vspace{.5cm}

\begin{theo}\label{theoa04}  If $\d$ is a fuzzy topology on $X$, then $\i(\d)$ is a topology on $X$ such that $\d \subset \om(\i(\d))$ with
equality if and only if $\d$ is affine invariant. Thus, if $\t = \i(\d)$, then $\d \subset \om(\t)$.
\end{theo}

{\bfseries Proof:} By definition of $\i(\d)$ every $f \in \d$ is a continuous map from $(X,\i(\d))$ to $I_r$. Hence, $\d \subset \om(\i(\d))$.
Thus, if $\t = \i(\d)$, then $\d \subset \om(\i(\d)) = \om(\t)$.

From Theorem \ref{theoa03} $\om(\t)$ is affine invariant for any topology $\t$ and so $\d = \om(\i(\d))$ implies that $\d$ is
affine invariant.

For the converse we prove a more precise result.

\begin{theo}\label{theoa05}If $\d$ is a fuzzy topology on $X$, then the
following conditions are equivalent.
\begin{enumerate}

\item[(a)]  $\d$ is affine invariant.

\item[(b)] $\d$ is laminated and if  $x \in U$ and $ U \in \i(\d)$, then there exists
$f \in \d $ such that $f(x) = 1$ and $f|(X \setminus U) = 0$.

\item[(c)] $\d$ is laminated and if $ U \in \i(\d)$, then $\chi_U \in \d$
where $\chi_U$ is the characteristic function of $U$.

\item[(d)] $\d = \om(\i(\d))$.

\item[(e)]  $\d$ is topological, i.e. $\d$ is induced from some topology $\t$ on $X$.
\end{enumerate}
\end{theo}

{\bfseries Proof:}  (a) $\Rightarrow$ (b):  If $U \in \i(\d)$ and $x \in U$, then by
definition of a subbase there exists a finite index set $K$, elements
$\{ g_i \in \d : i \in K \}$ and $\{ a_i \in I : i \in K \}$ such that
\begin{equation}\label{02}
x \ \in \ \bigcap_{i \in K} \ (g_i)_{(a_i} \ \subset \ U.
\end{equation}

By (a) condition (i'') holds for $\d$ and so, in particular, (i') holds as well, i.e. $\d$ is laminated.

Let $b_i = g_i(x)$ so that $0 \leq a_i < b_i \leq 1$. By (i'),(ii) and (iii)
$h_i =_{def} a_i \vee ( b_i \wedge g_i) \in \d$ and maps $X$ to the interval
$[a_i,b_i]$.  By (i'')
$f_i = (h_i  -  a_i)/(b_i - a_i)$ is in $\d$.
Furthermore, $f_i(x) = 1$ for all $i \in K$ and $\bigcap_{i \in K}  (f_i)_{(0}  \subset  U$.
Finally, (iii) implies that $f = \wedge_{i \in K} \ f_i  \ \in \ \d$. Clearly, $f(x) = 1$ and $f|(X \setminus U) = 0$.

(b) $\Rightarrow$ (c):  For each $x \in U$, (b) implies there exists $f_x \in \d$ such that $f_x(x) = 1$ and
$f_x|(X \setminus U) = 0$. By (ii) $\chi_U = \bigvee_{x \in U} f_x$ is in $\d$.

(c) $\Rightarrow$ (d): Assume that $f : X \to I$ is lsc with respect to $\i(\d)$. We must show that
$f \in \d$.

For each $c \in I$ the set $U_c = f_{(c}$ is open, i.e. lies in $\i(\d)$,
 because $f$ is lsc.  Hence, by (c),
the characteristic function $\chi_{U_c} \in \d$. Let $f_c = c \wedge \chi_{U_c}$ which lies in $\d$ by
(i') and (iii). Since $f$ maps to $I, \ f \geq f_c$.  Hence,
\begin{equation}\label{03}
f \ \geq \ \hat{f} \ =_{def} \ \bigvee_{c \in I} \ f_c.
\end{equation}

By (ii) $\hat{f} \in \d$.

If $f(x) = 0 $ then of course $ \hat{f}(x) \geq f(x)$.

Now assume that $f(x) > 0$ and let $\ep > 0$. Choose $c \in I$ such that $f(x) > c > f(x) - \ep$.
Since $f(x) > c, \ x \in U_c$ and so $f_c(x) = c$.  Hence, $\hat{f}(x) \geq f_c(x) > f(x) - \ep$.
As $\ep$ was arbitrary we have $\hat{f}(x) \geq f(x)$.

Thus, $f = \hat{f}$ and so $f \in \d$ as required.

(d) $\Rightarrow$ (e): Obvious.

(e) $\Rightarrow$ (a):  By Theorem \ref{theoa03} an induced fuzzy topology is affine invariant.

$\Box$ \vspace{.5cm}

{\bfseries Remark:} Notice that in the proof of (a) $\Rightarrow$ (b) we did not need the full strength
of the affine invariance assumption. We needed that $\d$ be laminated and if
 $f \in \d$ has image $f(X) \subset [a,b]$ with $0 \leq a < b \leq 1$,  then
$(f - a)/(b - a)$ is in $\d$.
\vspace{.5cm}

Following Martin \cite{M80} we define for a fuzzy topology $\d$ on $X$ the topology $\chi^*(\d) = \{ U  \subset X : \chi_U \in \d \}$ on $X$.
Since $U = (\chi_U)_{(0}$ it follows that $\chi^*(\d) \subset \i(\d)$. For a topology $\t$ we define the fuzzy topology
$\chi(\t) = \{ \chi_U : U \in \t \}$. Clearly, $\chi(\t) \subset \om(\t)$ and  $\t = \i(\chi(\t))$.

\begin{df}\label{defa07} We call a fuzzy topology $\d$ a \emph{weakly induced fuzzy topology} when it satisfies the
equivalent conditions
\begin{enumerate}
\item[(a)]$\chi^*(\d) = \i(\d)$.
\item[(b)]$\d \supset \chi(\i(\d))$.
\item[(c)] $\d \subset \om(\chi^*(\d))$.
\item[(d)] If $U = f_{(c}$ for $f \in \d$, then $\chi_U \in \d$.
\end{enumerate}
\end{df}

\vspace{.5cm}

Observe that condition (d) is just a restatement of conditions (a), (b) and (c).

Martin calls $\d$ an induced fuzzy topology  when $\d = \om(\chi^*(\d))$. We show that this definition
agrees with the one given in Definition \ref{defa02}.

\begin{prop}\label{propa06}If $\d$ is a fuzzy topology on $X$, then the
following conditions are equivalent.
\begin{enumerate}

\item[(a)] $\d = \om(\chi^*(\d))$.
\item[(b)] $\d$ is laminated and weakly induced.
\item[(c)] $\d = \om(\i(\d))$.
\item[(d)]  $\d$ is topological, i.e. $\d$ is induced from some topology $\t$ on $X$.
\item[(e)]  $\d$ is affine invariant.
\end{enumerate}
\end{prop}

{\bfseries Proof:} (a) $\Rightarrow$ (b): By Theorem \ref{theoa03} $\d = \om(\chi^*(\d))$ implies $\d$ is
laminated and $\chi^*(\d)= \i(\d)$. Hence, $\d$ is weakly induced.

(b)$\Leftrightarrow$ (c):  Conditions (b) and (c) are exactly conditions (c) and (d) of Theorem \ref{theoa05} and so the equivalence follows
from that theorem.

(c)$\Leftrightarrow$ (d)$\Leftrightarrow$ (e): These are restatements of parts of Theorem \ref{theoa05}.

(c) $\Rightarrow$ (a): We have seen that  (c) $\Rightarrow$ (b). Hence,  $\d = \om(\i(\d))$ implies that $\chi^*(\d) = \i(\d)$
and so that $\d = \om(\chi^*(\d))$.

$\Box$ \vspace{.5cm}

\begin{theo}\label{theoa09} Let $\t$ be a topology on $X$.
\begin{enumerate}

\item[(a)] The affine invariant fuzzy topology $\om(\t)$  is generated by \\
$\{ k : k \in I \} \cup \chi(\t) = \{ k : k \in I \} \cup \{ \chi_U : U \in \t \}.$

\item[(b)] If $\t$ is completely regular, then $\om(\t)$  is generated by the continuous functions in $I^X$.
\end{enumerate}
\end{theo}

{\bfseries Proof:} (a): Let $\d$ be the fuzzy topology generated by  $\{ k : k \in I \} \cup \chi(\t)$. Since
the generators are lsc, it follows that $\d \subset \om(\t)$. Hence, $\i(\d) \subset \i(\om(\t)) = \t$. Clearly, $U \in \i(\d)$
for all $U \in \t$ and so $\i(\d) = \t$. Hence, $\d$ satisfies condition (c) of Theorem \ref{theoa05}. From the theorem it follows
that $\d = \om(\i(\d)) = \om(\t)$.

(b): Let $\d$ be the fuzzy topology generated by the continuous functions in $I^X$. Again, $\d \subset \om(\t)$ and $\i(\d) \subset \i(\om(\t)) = \t$.
Because $\t$ is completely regular, the collection $\{ f_{(c} \}$ with $c \in I$ and $f : X \to I$ continuous form a basis for $\t$ and
so $\t \subset \i(\d)$. Thus, $\t = \i(\d)$. By complete regularity $\d$ satisfies condition (b) of Theorem \ref{theoa05}. From the theorem it follows
that $\d = \om(\i(\d)) = \om(\t)$.

$\Box$ \vspace{.5cm}

If $H : X_1 \to X_2$ then for any  $L \subset \R$ we define $H^* : L^{X_2} \to L^{X_1}$ by $H^*(f) = f \circ H$.
If $(X_1,\d_1)$ and $(X_2,\d_2)$ are fuzzy topological spaces then $H $ is defined to be \emph{fuzzy continuous}
when $H^*(\d_2) \subset \d_1$. We thus obtain the category of fuzzy topological spaces $FTOP$ and the full subcategories of
laminated fuzzy topological spaces $FLTOP$ and affine invariant fuzzy topological spaces $FATOP$.

A map  $H : (X_1,\t_1) \to (X_2,\t_2)$ between topological spaces is a \emph{quotient map} when $U \in \t_2$ if and only if $H^{-1}(U) \in \t_1$.
A map  $H : (X_1,\d_1) \to (X_2,\d_2)$ between fuzzy topological spaces is a \emph{fuzzy quotient map} when $f \in \t_2$ if and only if $H^*(f) \in \d_1$.

We review the proof of the following well-known result.

\begin{theo}\label{theoa10}  For $ i = 1,2$ let $\t_i$ and $\d_i$ be a topology and a fuzzy topology on a set $X_i$.
Let $H : X_1 \to X_2$ be a set map.

(a) If $H : (X_1,\d_1) \to (X_2,\d_2)$ is fuzzy continuous, then $H : (X_1,\i(\d_1)) \to (X_2,\i(\d_2))$
is continuous.

(b) $H : (X_1,\t_1) \to (X_2,\t_2)$ is continuous if and only if $H : (X_1,\om(\t_1)) \to (X_2,\om(\t_2))$ is fuzzy
continuous. $H : (X_1,\t_1) \to (X_2,\t_2)$ is a quotient map if and only if $H : (X_1,\om(\t_1)) \to (X_2,\om(\t_2))$ is fuzzy
quotient map.

(c) If $\t_1$ is affine invariant, then $H : (X_1,\d_1) \to (X_2,\d_2)$ is fuzzy continuous if and only if
$H : (X_1,\i(\d_1)) \to (X_2,\i(\d_2))$.
\end{theo}

{\bfseries Proof:}  If $f \in I^{X_2}$ and $c \in I$, then
\begin{equation}\label{04}
H^*(f)_{(c} \quad = \quad H^{-1}(f_{(c}).
\end{equation}

(a): If $H$ is fuzzy continuous, then $f \in \d_2$ implies $H^*(f) \in \d_1$ and so equation (\ref{04}) says that
$H$ pulls the elements of the defining subbase for $\i(\d_2)$ back into $\i(\d_1)$.  Hence, $H$ is
continuous.

(b): If $H$ is continuous and $f$ is lsc, then $H^*(f) = f \circ H$ is lsc by equation (\ref{04}) and if $H$ is a quotient map,
then the converse is true as well. Hence,
$H^*$ maps $\om(\t_2)$ to $\om(\t_1)$ and so $H$ is fuzzy continuous and $H$ is a fuzzy quotient map if it is a quotient map. If $H$ is fuzzy continuous,
then by (a) $H : (X_1,\i(\om(\t_1))) \to (X_2,\i(\om(\t_2)))$ is continuous and so, $H$ is continuous by
Theorem \ref{theoa03}. Now assume that $H$ is a fuzzy quotient map and $U \subset X_2$ with $H^{-1}(U) \in \t_1$. Thus,
$H^*(\chi_U) = \chi_{H^{-1}(U)} \in \om(\t_1)$. Because $H$ is a fuzzy quotient map, $\chi_U \in \om(\t_2)$, i.e. it is lsc with respect to $\t_2$.
Hence, $U \in \t_2$.

(c): If $H$ is  continuous, then by (b) $H : (X_1,\om(\i(\d_1))) \to (X_2,\om(\i(\d_2)))$ is fuzzy
continuous.  If $\d_1$ is affine invariant, then by Theorem \ref{theoa04}
$\d_2 \subset \om(\i(\d_2))$ and  $\d_1 = \om(\i(\d_1))$. Thus, $H$ is fuzzy continuous.

$\Box$  \vspace{.5cm}

We thus have functors $\tilde \i : FTOP \to TOP$ and $\tilde \om : TOP \to FTOP$ with $\tilde \i \circ \tilde \om$ the
identity on $TOP$ and $\tilde \om \circ \tilde i$ a retraction from $FTOP$ onto $FATOP$.

\begin{prop}\label{propa10a} For the unit interval $I_r$ equipped with the lower topology $\t_r$,
let $\d_r = \om(\t_r)$ be the induced fuzzy topology.
Let $(X,\d)$ be a fuzzy topological space. If $f : (X,\d) \to (I_r,\d_r)$ is fuzzy continuous then $f \in \d$. The fuzzy topology $\d$
is affine invariant if and only if  $f : (X,\d) \to (I_r,\d_r)$ is fuzzy continuous for all $f \in \d$. \end{prop}

{\bfseries Proof:} The identity map $1_I$ is a continuous function on $I_r$ and so $1_I \in \d_r$.
Hence, if $f :(X,\d) \to (I_r,\d_r)$ is fuzzy continuous, then $f^*(1_I) = f \in \d$.

If $m, k \in \R$ with $m > 0$, let
\begin{equation}\label{04a}
\phi_{m,k}(t) \ = \ \max(0,\min(1, m \cdot t + k)).
\end{equation}
Clearly, $\phi_{m,k} : I \to I$ is continuous and non-decreasing and so $\phi_{m,k} : (I_r,\t_r) \to (I_r,\t_r)$ is continuous.
Hence, $\phi_{m,k} \in \d_r$.
So if $f :(X,\d) \to (I_r,\d_r)$ is fuzzy continuous then $f^*(\phi_{m,k}) = \phi_{m,k} \circ f \in \d$. Thus, if
every $f \in \d$ is fuzzy continuous then $\d$ is affine invariant.

If $\d$ is affine invariant, then with $\t = \i(\d), \d = \om(\t)$. So if $f \in \d$ then $f$ is lsc with respect to $\t$, i.e.
$f : (X,\t) \to (I_r,\t_r)$ is continuous and so by Theorem \ref{theoa10} (b),  $f : (X,\d) \to (I_r,\d_r)$  is fuzzy continuous.

$\Box$  \vspace{.5cm}

{\bfseries Remark:} The identity from $\R$ to $\R_r$ is continuous and so is fuzzy continuous. Hence, if $f : (X,\d) \to (I,\d_I)$
is fuzzy continuous with $\d_I$ the fuzzy topology induced from the usual topology $\t_I$ on $I$
then $f : (X,\d) \to (I_r,\d_r)$ is fuzzy continuous.
\vspace{.5cm}

Let $Y \subset X$. If $\t$ is a topology on $X$, then the \emph{relative topology} $\t|Y$ on $Y$ is $\{ U \cap Y : U \in \t \}$.
For $f \in I^X$, let $f|Y$ denote the restriction of $f$ to $Y$ so that $f|Y \in I^Y$.   If $\d$ is a fuzzy topology on $X$,
 then the \emph{relative fuzzy topology} $\d|Y$ on $Y$ is $\{ f|Y : f \in \d \}$.

\begin{prop}\label{propa10b} Let $\t$ and $\d$ be a topology and a fuzzy topology on $X$ and let $Y \subset X$.
\begin{equation}\label{04b}
\i(\d|Y) \ = \ (\i(\d))|Y \qquad \text{and} \qquad \om(\t|Y) \ = \ (\om(\t))|Y
\end{equation}
If $\d$ is laminated or affine invariant then $\d|Y$ satisfies the corresponding property.
\end{prop}

{\bfseries Proof:} It is clear that $\d|Y$ is a fuzzy topology which is laminated if $\d$ is. Let $f \in \d$ and
$m, k \in \R$ so that with $\phi(t) = m \cdot t + k$, $\phi(f(y)) \in I$ for all $y \in Y$. Since $\phi$ is continuous and
increasing $\phi^{-1}(I)$ is a closed interval in $\R$.  If $b = \sup f(Y), a = \inf f(Y)$, then $\phi([a,b]) \subset I$.
If $\bar f = a \vee (b \wedge f)$ then  $\bar f|Y = f|Y$ and $\phi(f(x)) \in I$ for all $x \in X$. If $\d$ is affine invariant,
then $\bar f \in \d$ and $\phi \circ \bar f \in \d$ and so that $ \phi \circ (f|Y) = (\phi \circ \bar f)|Y \in \d|Y$.  That is, $\d$ is affine invariant.

Since $f_{(c} \cap Y = (f|Y)_{(c}$ and $A \mapsto A \cap Y$ commutes with unions and intersections, it follows that
$\i(\d|Y)  =  (\i(\d))|Y$.

Since $(\om(\t))|Y$ is affine invariant, Theorems \ref{theoa04} and \ref{theoa03} imply
\begin{equation}\label{04c}
 (\om(\t))|Y \ = \ \om(\i[(\om(\t))|Y]) \ = \ \om([\i(\om(\t))]|Y) \ = \ \om(\t|Y).
 \end{equation}

$\Box$  \vspace{.5cm}

{\bfseries Remark:}  In particular, this shows that any bounded lsc function on $Y$ extends to an lsc function on $X$.
\vspace{.5cm}

For $ \{ X_i : i \in K \}$  an indexed family of sets let $\prod_{i \in K} \ X_i$ denote the product with projections
$\pi_j :  \prod_{i \in K} \ X_j \to X_j$ for all $j \in K$.
If $\t_i$ is a topology on $X_i$ then we denote by
$\prod_{i \in K} \ \t_i$ the product topology on $\prod_{i \in K} \ X_i$.  If $\d_i$ is a fuzzy topology on $X_i$ then
we denote by $\prod_{i \in K} \ \d_i $ the fuzzy topology on $\prod_{i \in K} \ X_i$
generated by $\bigcup_{i \in K} \ \pi_i^*(\d_i)$.  This is the coarsest fuzzy topology on the product so that
the projection maps $\pi_i$ are fuzzy continuous.

\begin{prop}\label{propa11} If $\{ (X_i, \d_i) : i \in K \}$ is an indexed family of fuzzy topological spaces, then we have
\begin{equation}\label{05}
\i(\prod_{i \in K} \ \d_i) \ = \  \prod_{i \in K} \ \i(\d_i).
\end{equation}

If  $\d_i$ is laminated for some $i \in K$ then $\prod_{i \in K} \ \d_i $ is laminated.

If  $\d_i$  is affine invariant for all $i \in K$ then $\prod_{i \in K} \ \d_i$ is affine invariant.
\end{prop}

{\bfseries Proof:}  If $f \in \d_i$ and $c \in I$ then
$\pi_i^{-1}(f_{(c}) = (\pi_i^*(f))_{(c} \in \i(\prod_{i \in K} \ \d_i)$. Hence,
$\prod_{i \in K} \ \i(\d_i) \subset \i(\prod_{i \in K} \ \d_i)$.

On the other hand, the set of lsc functions $\om(\prod_{i \in K} \ \i(\d_i))$ is a fuzzy topology
which clearly contains $\bigcup_{i \in K} \ \pi_i^*(\d_i)$ and so contains $\prod_{i \in K} \ \d_i$.
So Theorem \ref{theoa03} implies
\begin{equation}\label{06}
\prod_{i \in K} \ \i(\d_i) \ = \ \i(\om(\prod_{i \in K} \ \i(\d_i))) \ \supset \ \i(\prod_{i \in K} \ \d_i).
\end{equation}

If $\d_j$ is laminated then the constant functions from $\prod_{i \in K} \ X_i$ to $I$
are contained in $\pi_j^*(\d_j)$.

If all the $\d_i$'s are affine invariant then they all laminated
 and so we can apply Theorem \ref{theoa05} to $\prod_{i \in K} \ \d_i$
and check that the product fuzzy topology satisfies condition (b) of Theorem \ref{theoa05}.

If $x \in U$ and $U$ is open in the product topology then there is a
finite subset $\hat{K}$ of $K$ and open sets $\{ U_i : i \in \hat{K} \}$
so that $x \in \hat{U} =_{def} \bigcap_{i \in \hat{K}} \pi_i^{-1}(U_i) \subset U$.
By Theorem \ref{theoa05} condition (c), each $\chi_{U_i} \in \d_i$.
Hence,
\begin{equation}\label{07}
f \ =_{def} \  \chi_{\hat{U}} \ = \ \bigwedge_{i \in \hat{K}} \ \pi_i^*(\chi_{U_i}) \ \in \prod_{i \in K} \ \d_i.
\end{equation}
Clearly, $f(x) = 1$ and $f|(X \setminus U) = 0$.

$\Box$ \vspace{.5cm}

The coproduct $\coprod_{i \in K} \ X_i$ is defined to be the disjoint union of
the indexed sets $X_i$ so that each $X_i$ can be
regarded as a subset of $\coprod_{i \in K} \ X_i$.  If $\t_i$ is a topology on
$X_i$ then the topology $\coprod_{i \in K} \ \t_i$ is the collection of subsets $U$
such that $U \cap X_i \in \t_i$ for all $i \in K$.
If $\d_i$ is a fuzzy topology on $X_i$ then the fuzzy topology $\coprod_{i \in K} \ \d_i $
is the set of functions $f$ such that
$f|X_i \in \d_i$ for all $i \in K$.  In each case the required axioms are easy to check.
We leave the easy proof of the following analogue of
Proposition \ref{propa11} to the reader.

\begin{prop}\label{propa12}  If $\{ (X_i, \d_i) : i \in K \}$ is an indexed family of fuzzy topological spaces, then
\begin{equation}\label{08}
\i(\coprod_{i \in K} \ \d_i) \ = \ \coprod_{i \in K} \ \i(\d_i).
\end{equation}

If  $\d_i$  is laminated for all $i \in K$ or is affine invariant for all
$i \in K$ then $\coprod_{i \in K} \ \d_i$ satisfies the corresponding property.
\end{prop}

$\Box$  \vspace{1cm}

\begin{cor}\label{cora13}  If $\{ (X_i, \t_i) : i \in K \}$ is an indexed family of  topological spaces and $\d_i = \om(\t_i)$ is
the induced fuzzy topology  on $X_i$, then
\begin{align}\label{09}
\begin{split}
\prod_{i \in K} \d_i \ = \ &\om(\prod_{i \in K} \t_i), \\
\coprod_{i \in K} \d_i \ = \ &\om(\coprod_{i \in K} \t_i),
\end{split}
\end{align}
are the induced   fuzzy topologies on the product and coproduct.
\end{cor}

 {\bfseries Proof:}  By  Theorem \ref{theoa03}  $\t_i = \i(\d_i)$. By Proposition \ref{propa11} the product
fuzzy topology is affine invariant. Hence, from (\ref{05}) and Theorem \ref{theoa04} we have that
\begin{align}\label{10}
\begin{split}
\prod_{i \in K} \ \d_i \ = \ &\om(\i(\prod_{i \in K} \ \d_i)) \ = \\
\om(\prod_{i \in K} \ \i(\d_i)) \ &= \  \om(\prod_{i \in K} \ \t_i),
\end{split}
\end{align}
 as required.

 The proof for the coproduct is similar using Proposition \ref{propa12} and (\ref{08}).

 $\Box$ \vspace{.5cm}

We conclude by illustrating some pathologies that can occur without lamination or affine invariance.
\vspace{.5cm}

{\bfseries Examples:}  (A) For any topological space $(X, \t)$ the set of
characteristic functions of the open sets is the fuzzy topology $\chi(\t)$
which is not laminated, i.e. does not satisfy (i').  In fact,
Let $L$ be a proper  subset of $I$ which contains $0$ and $1$. Call $L$ \emph{sup-closed} if $\sup A \in L$ whenever
$A$ is a nonempty subset of $L$.  If $L$ is closed, then, of course, it is sup-closed. For any fuzzy topology $\d$ the
set of $k \in I$ such that the constant $k \in \d$ is sup-closed.
If $L$ is a proper, sup-closed  subset of $I$ which contains $0$ and $1$.
then $\om_L(\t) =_{def} L^X \ \cap \ \om(\t)$, the set of lsc
functions with image in $L$, is a fuzzy topology which is not laminated. Only the constants $k$ with $k \in L$
are in $\d$. Clearly, each
$\d = \om_L(\t)$ satisfies $\t = \i(\d)$ and is weakly induced, but is not, of course, induced. We sketch the proof of
the following extension of Proposition \ref{propa06} and Theorem \ref{theoa09} (a).

\begin{prop}\label{propa14} (a) Let $\d$ be a fuzzy topology and $L$ be the set of constants $k \in \d$.
 If $\d$ is weakly induced, then $\d \cap L^X = \om_L(\i(\d))$.

 (b) If $\t$ is a topology and $L$ is a sup-closed subset which contains $0$ and $1$ then $\om_L(\t)$ is generated
 by $L \cup \chi(\t)$.
 \end{prop}

{\bfseries Proof:} (a): Since $\d \subset \om(\i(\d))$ it is clear that $\d \cap L^X \subset \om_L(\i(\d))$.
Now adapt the Theorem \ref{theoa05} proof of (c) $\Rightarrow$ (d). Let $f \in L^X$ be lsc
with respect to $\i(\d)$. Call $c \in L$ a right endpoint if
$L \cap (c - \ep,c) = \emptyset$ for some $\ep > 0$, i.e. $c$ is a right endpoint of one of the maximum open subintervals of
$I \setminus L$.  In that case, $f_{[c} \in \i(\d)$. For $c$ a right endpoint in $L$
let $f_c = c \wedge \chi_{f_{[c}}$.  For the remaining $c \in L$ let $f_c = c \wedge \chi_{f_{(c}}$. As before,
$f = \bigvee_{c \in L} \ f_c \in \d$.

(b): The proof is the obvious analogue of that of Theorem \ref{theoa09} (a) using part (a) here in place of Theorem \ref{theoa05}.

 $\Box$ \vspace{.5cm}

 If $Y$ is a proper open subset of $X$ then we can extend any lsc function $f : Y \to I$ to an lsc function on
 $X$ by letting $f(x) = 0$ for $x \in X \setminus Y$.  Thus, we can include $\om(\t|Y) \subset \om(\t)$. If
 $\d$ is the fuzzy topology on $X$ generated by $\om(\t|Y) \cup \chi(\t)$ then $\d = \{ f \vee \chi_U :
 f \in \om(\t|Y), U \in \t \}$. Thus, $\d$ is a weakly induced fuzzy topology with $\t = \i(\d)$ and with
 no constants other than $0$ and $1$.

 (B) We say that a subinterval $J$ of $I$ is \emph{non-trivial} when it has positive length.
We say that two  subintervals $J, L$ of $I$ are
\emph{non-overlapping} when their interiors are disjoint. When
$J$ and $L$ are non-overlapping then either $J \leq L$, meaning
$s \leq t$ for all $s \in J$ and $t \in L$, or else $L \leq J$.
Notice that if $J \leq L$ then for all $f \in J^X, g \in L^X$
\begin{equation}\label{11}
f \wedge g \ = \ f \qquad \mbox{and} \qquad f \vee g \ = \ g.
\end{equation}

Let $\J = \{ J_1, J_2,...  \}$ be a finite or infinite sequence of
nontrivial, closed subintervals of $I$ with any two non-overlapping. For a set $X$
let
\begin{equation}\label{12}
\d_{\J} \ =_{def} \  I \ \cup \ \bigcup_{i = 1,2,...} \ J_i^X \  \subset \  I^X.
\end{equation}
where the $I$ in the union denotes the set of constants in $I$. From (\ref{11}) it follows that
 $\d_{\J}$ is a fuzzy topology.

For example, if $\J = \{ I \}$ then $\d_{\J}(X) = I^X$.

For a topological space $(X,\t)$ let
\begin{equation}\label{13}
\om_{\J}(\t) \quad =_{def} \quad \d_{\J} \cap \om(\t).
\end{equation}

 It is clear that any $\om_{\J}(\t)$ is
a laminated fuzzy topology.  Furthermore, if $a, b \in J_1$ with $a < b$ then
for any $U \in \t, \ f = a \vee(b \wedge \chi_U) $
has image in $J_1$ and so lies in $\om_{\J}(\t)$. Also, $U = f_{(a}$. Thus, $\t = \i(\om_{\J}(\t))$.

If $\J $ is not $\{ I \}$, then $\om_{\J}(\t)$ is a proper subset of $\om(\t)$.  In fact, for $U \in \t$,
$\chi_U \not\in \om_{\J}(\t)$ unless $U = X$ or $\emptyset$ in which case $\chi_U$ is constant.
Thus, $\chi^*(\om_{\J}(\t)) = \{ \emptyset, X \}$.

Furthermore, if $(X, \t) = \coprod_{k \in K} \ (X_k, \t_k)$ then we can choose a separate
sequence $\J_k$ for each $k \in K$ and define a laminated fuzzy topology  on $X$ by using
$\coprod_{k \in K} \ \om_{\J_k}(\t_k)$.

Thus, we
obtain many examples of laminated fuzzy topologies  each of which generates
 $\t$ but which are not affine invariant.

(C) The real pathologies occur with products.  Suppose $(X_1, \t_1)$ and $(X_2, \t_2)$ are topological spaces
and $(X, \t) = (X_1 \times X_2, \t_1 \times \t_2)$. On $X_1$ let $\d_1 = \om_{\{[0,1/2] \} }(\t_1)$ and on
$X_2$ let $\d_2 = \om_{\{ [1/2,1] \} }(\t_2)$. Using equation (\ref{11}) again it is easy to check that the union
\begin{equation}\label{14}
\d \ =_{def} \  \pi_1^* \d_1  \ \cup \ \pi_2^* \d_2
\end{equation}
 already satisfies conditions (i'),(ii) and (iii) of Definition \ref{defa01}
 and so is the product fuzzy topology.
 By Proposition \ref{propa11}, $\i(\d)$ is the product topology $\t$.  However, each element of $\d$ is a function depending on
 either the first or the second variable alone.  That is, each fuzzy set is horizontal or vertical.

(D) Finally, let $(X,\t)$ be a second countable topological space with $\B$ any countable subbase for $\t$.
Let $\J$ be a countably infinite set of nontrivial, closed subintervals of $I$ with any two non-overlapping,
e.g the closures of the components of the complement of the Cantor set. Let $\rho : \B \to \J$ be a
bijection.  For $U \in \B $ we will write $ \rho(U) = [a(U),b(U)]  $ so that $0 \leq a(U) < b(U) \leq 1$. Let
%\begin{equation}
%\G_{\rho} \quad =_{def} \quad \{  a(U) \vee (b(U) \wedge \chi_U) : U \in \B \ \} \cup \{ 0,1 \}. \hspace{1cm}
%\end{equation}
%From equation (1.10) again it follows that the countable set $\G_{\rho}$ is a fuzzy topology.  Furthermore,
%$\T = \T(G_{\rho})$.  Of course, $\G_{\rho}$ is not laminated.  However, we can define
\begin{equation}\label{15}
\d_{\rho} \  =_{def} \  \{ c \vee (d \wedge \chi_U) : U \in \B \ \& \ a(U) \leq c < d \leq b(U) \} \cup I.
\end{equation}

Now $\d_{\rho}$  is a laminated fuzzy topology with $ \t = \i(\d_{\rho})$ but nonetheless
every element of $\d_{\rho}$ is just an affine adjustment of the characteristic function of a set in $\B$.
Hence,
\begin{equation}\label{16}
\{ f_{(c} : f \in \d_{\rho}, \ c \in I \} \ = \ \B \cup \{ X, \emptyset \}.
\end{equation}

If, in addition, $(X,\t)$ is connected then the constants are the only elements of $\d_{\rho}$ which are
continuous. Let  $I$ be equipped with the usual topology $\t_I$
and the induced fuzzy topology $\d_I$. By Proposition \ref{propa10b} and the Remark thereafter,
if $H :(X,\d_{\rho}) \to (I,\d_I)$ is fuzzy continuous, then $H \in \d_{\rho}$. Furthermore, Theorem \ref{theoa10}(a)
implies that $H :(X,\t) \to (I,\t_I)$ is continuous.
Hence, the constant functions are the only fuzzy continuous maps from
$(X,\d_{\rho})$ to $(I,\d_I)$ when $(X,\t)$ is connected.

(E) Assume that $(X_1, \d_1)$ and $(X_2, \d_2)$ are fuzzy topological spaces and that
$\d_2$ is laminated.
Every constant  $k \in \d_2$.
If $H : X_1 \to X_2$ is any fuzzy continuous function
then $k = H^*k  \in \d_1$ and so $\d_1$ is laminated
 as well. Thus, if $\d_1$ is not laminated, then there is no fuzzy continuous function from $(X_1,\d_1)$
to $(X_2,\d_2)$.

 \vspace{1cm}

\section{Compactness}

From now on for topological spaces $X, Y$ etc, we suppress explicit labeling of the topology,
and equip each with the induced fuzzy topology, i.e. for $X$ with topology $\t$ we use the fuzzy topology
 $\d = \om(\t)$, the set of lsc functions in $I^X$.  We will call these the fuzzy open  sets.  Thus,   $f$ is a fuzzy open
 set if and only if $f_{(c}$ is  open for every $c \in I$.  On the other hand, $f$ is fuzzy closed when its complement $1 - f$ is fuzzy
 open and so when
 $f$ is usc.  Thus, $f$ is a fuzzy closed set when
 $f_{[c}$ is closed for every $c \in I$.

 We turn now to compactness.  Since we are not - yet - assuming the spaces are Hausdorff, a compact set
 need not be closed though a closed subset of a compact set is compact.

Among the possible definitions for fuzzy compactness, we follow Weiss \cite{W75}

 \begin{df}\label{defc01}  Let $X$ be a topological space with  the induced fuzzy
 topology  and let $f$ be a fuzzy set, i.e. $f \in I^X$.
\begin{enumerate}
 \item[(a)] We say that $f$ is \emph{fuzzy compact} when $f_{[c}$ is compact for all $ c \in (0,1]$.

  \item[(b)] We say that $f$ satisfies \emph{ Condition L} if whenever $\{ g_i : i \in K \}$ is an indexed
 family of fuzzy open sets such that $\bigvee_{i \in K} \ g_i  \geq f$ and $\ep > 0$  there
 exists a finite  $\hat{K} \subset K$ such that $\bigvee_{i \in \hat{K}} \ g_i  \geq f - \ep$.
 \end{enumerate}
 \end{df}
 \vspace{.5cm}

 In the definition of fuzzy compactness we do not assume that $ f_{[0}$, which is all of
 $X$, is compact.  In particular, if $f$ is fuzzy closed and has compact support, i.e. the closure of $f_{(0}$
 is compact, then $f$ is  fuzzy compact. Thus, the value $0$ plays a special role and fuzzy compactness need not
 be preserved by affine adjustments.

 Condition L was introduced in Lowen \cite{L76} as a definition of  compactness in the context of
 general laminated fuzzy topologies.

 \begin{lem}\label{lemc02} Let $f $ be a fuzzy subset of a topological space $X$. If for every
  $c \in (0,1)$ the set $f_{[c}$ is compact, then $f$ satisfies
 condition L. Furthermore, if $g$ is fuzzy closed and $g \leq f$, then $g$ is fuzzy compact.
 \end{lem}

 {\bfseries Proof:}  To verify Condition L we begin with a family $\{ g_i : i \in K \}$ of fuzzy open
 sets with $\bigvee_{i \in K} \ g_i \ \geq \ f$ and an arbitrary $\ep > 0$.

 Let $1 = c_0 > c_1 > ... > c_{n} = 0$ be a decreasing sequence
 in $I$ with $ c_{k-1} - c_k < \ep/2$  for $k = 1,.., n$.

The set $\{ (g_i)_{(c_k}  : i \in K \}$ is an open cover of the compact set $f_{[c_{k-1}}$
for  $k = 2,..., n$.  Choose $K_k \subset K$ finite so the  $\{ (g_i)_{(c_k}  : i \in K_k \}$
is a subcover.

Let $\hat{K} = \bigcup_{k = 2,..,n} K_k$.

If for some $k = 2,...,n \ f(x) \in [c_{k-2},c_{k-1}] $ then $g_i(x) > c_k \geq f(x) - \ep $ for some $i \in K_k$.

If $f(x) \in (c_{n-1},c_n] $ then $f(x) - \ep < 0 \leq g_i(x)$ for any $i$.

Thus, as required  $\bigvee_{i \in \hat{K}} \ g_i  \geq f - \ep$.

Now if $g \leq f$ and $g$ is fuzzy closed, let $c \in (0,1]$. Choose $a$ so that $0 < a < c$. Observe that
$g_{[c}$ is a closed subset of $f_{[a}$ which is compact by assumption (since $a < 1$).  Hence, $g_{[c}$ is compact.

$\Box$ \vspace{.5cm}

\begin{theo}\label{theoc03}  Let $X$ be a topological space and let $f $ be a fuzzy subset of $X$.
\begin{enumerate}
\item[(a)] If $f$ is fuzzy compact then $f$ satisfies Condition L.

\item[(b)] If $f$ is fuzzy closed and satisfies Condition L then $f$ is fuzzy compact.
\end{enumerate}
\end{theo}

{\bfseries Proof:} (a):  This follows from Lemma \ref{lemc02}.

(b): Assume that $f$ is fuzzy closed and satisfies Condition L.
Let $c \in (0,1]$ and choose $\ep \in (0,c)$.

Assume $\{ U_i : i \in K \}$ is an open cover of $f_{[c}$.  Let $V = X \setminus f_{[c}$
which is open because $f$ is usc.  Apply Condition L to  $\{ \chi_V \} \cup \{ \chi_{U_i} : i \in K \}$ and $f$.
We obtain the finite set $\hat{K} \subset \{ V \} \cup K$. If $x \in f_{[c}$ then $\chi_V(x) = 0 <  c - \ep \leq f(x) - \ep$. Hence, there
exists $i \in \hat K \setminus \{ V \}$ such that $ \chi_{U_i}(x) > f(x) - \ep  \geq c - \ep > 0$. That is, $x \in U_i$. Hence,
 $\{ U_i : i \in  \hat{K} \} $ is a finite subcover of $f_{[c}$.

$\Box$ \vspace{.5cm}

In particular we have the following, due to Lowen, \cite{L76}.

\begin{cor}\label{corc04} For a topological space $X$ the following are equivalent.
 \begin{itemize}
 \item[(i)] X is compact.
 \item[(ii)] The constant function $1$ is fuzzy compact.
 \item[(iii)] Every fuzzy closed set is fuzzy compact.
 \item[(iv)] There exists $k \in (0,1]$ such that the
 constant function $k$ satisfies Condition L.
 \end{itemize}
 \end{cor}

 {\bfseries Proof:}  (i) $\Leftrightarrow$ (ii): Obvious

 (ii) $\Rightarrow$ (iii): Any fuzzy closed set $f$ satisfies $f \leq 1$ and
 so $f$ is fuzzy compact when $1$ is by Lemma \ref{lemc02}.

 (iii) $\Rightarrow$ (iv): Every constant function $k$ is  fuzzy closed and
 so by (iii) it is fuzzy compact.  Condition L follows from Theorem \ref{theoc03} (a).

 (iv) $\Rightarrow$ (i): If $k$ satisfies Condition L then by Theorem \ref{theoc03} (b) it
 is fuzzy compact and $X = k_{[k}$ is compact.

 $\Box$ \vspace{.5cm}

\begin{cor}\label{corc05}  Assume $X$ is a Hausdorff topological space and let $f \in I^X$.

If $f$ is fuzzy compact then
it is fuzzy closed. Furthermore, $f$ is fuzzy compact if and only if it satisfies Condition L.

If $X$ is compact Hausdorff then $f$ is fuzzy compact if and only if it is fuzzy closed.
\end{cor}

{\bfseries Proof}:  If $f$ is fuzzy compact then for each $c > 0$, $f_{[c}$ is compact and
so is closed because $X$ is Hausdorff.  $f_{[0} = X$ is always closed.  Hence, $f$ is
fuzzy closed. In addition, Theorem \ref{theoc03}(a) implies that $f$ satisfies condition L.  For the converse
it suffices by Theorem \ref{theoc03}(b) to show that if $f$ satisfies condition L then it is fuzzy closed.
We prove the contrapositive,
adapting the usual proof that a compact subset of a Hausdorff space is closed.

Suppose that $f(x) = a < c$ but that $x$ is in the closure $A$ of $f_{[c}$.
For each $z \in A' =_{def} A \setminus \{ x \}$ choose disjoint open sets $U_z, W_z$ with $z \in U_z$ and
$x \in W_z$.  Let $V$ be the open set $X \setminus A$. Define :
\begin{equation}\label{17}
\begin{split}
\G_{in} \  =_{def} \  \{ \chi_{U_z} : z \in A' \ \} \hspace{4cm}\\
 g_{in} = \bigvee \{ g \in \G_{in} \}   \quad \mbox{and} \quad g_{out} = a \vee \chi_V. \quad \hspace{2cm}
\end{split}
\end{equation}

Choose $\ep < c - a$. Observe that $g_{out}$ and the elements of $\G_{in}$ are
fuzzy open.  Furthermore, if $z \not= x$ then $(g_{out} \vee g_{in})(z) = 1$ and so $ g_{out} \vee g_{in} \geq f$.
  On the other hand for every
$z \in f_{[c}$, we have $f(z) - \ep \geq c - \ep > a = g_{out}(z) $.  So if
Condition L were true of $f$  then applying it to $\G_{in} \cup \{ g_{out} \}$ would yield
 a finite subset $\hat{A} \subset A'$ such that
$\{ U_z : z \in \hat{A} \}$ covers $f_{[c}$.  But for any finite
$\hat{A} \subset A', \ \bigcup \{ U_z : z \in \hat{A} \}$ is disjoint
from $\bigcap \{ W_z : z \in \hat{A} \}$ which is a neighborhood of $x$ and so contains points of
$f_{[c}$. Hence,  $\{ U_z : z \in \hat{A} \}$ does not cover $f_{[c}$ and so  Condition L fails.

Finally, if $X$ is compact then any fuzzy closed set is fuzzy compact by Corollary \ref{corc04}.

$\Box$ \vspace{.5cm}

The fuzzy Tychonoff Theorem now follows from the usual Tychonoff Theorem.

\begin{theo}\label{theoc06} Let $\{ X_i : i \in K \}$ be an indexed family of topological spaces with
product space $X = \prod_{i \in K} \ X_i$. Let $f_i$ be a fuzzy subset of $X_i$ for each $i \in K$
and define $f = \bigwedge_{i \in K} \ f_i \circ \pi_i$. If each $f_i$ is  fuzzy compact, then
$f$ is fuzzy compact.
\end{theo}

{\bfseries Proof:}  For each $c \in (0,1]$
\begin{equation}\label{18}
f_{[c} \quad = \quad \prod_{i \in K} \ (f_i)_{[c}.  \hspace{3cm}
\end{equation}
That is, for $x \in X \ f(x) \ = \ inf_{i \in K} \ f_i(x_i) \  \geq \ c$  if and only if $f_i(x_i) \geq c$ for all $i \in K$.
As the product of compact spaces, $f_{[c}$ is compact.

$\Box$ \vspace{.5cm}

Clearly a fuzzy closed set with compact support is fuzzy compact.  For $X = \R$
with the usual topology and induced fuzzy topology let $f(t) = exp(-t^2)$. The element $f$ is fuzzy compact
but does not have compact support.

For a locally compact, Hausdorff space $X$ let $X^* = X \cup \{ p \}$ denote the one-point compactification
so that $U \subset X^*$ is open if and only if $U \cap X$ is open and, in addition, $p \not\in U$ or
$X \setminus U$ is compact.  $X^*$ is a compact Hausdorff space.
  For $f \in I^X$ let $f^* \in I^{X^*}$ be the extension of $f$ with $f^*(p) = 0$.  The following
  is then obvious.

\begin{prop}\label{propc07}  For a locally compact, Hausdorff space $X$ with one-point compactification $X^*$ a fuzzy set
$f \in I^X$ is a fuzzy compact set in $X$ if and only if the  extension $f^* \in I^{X^*}$ is a fuzzy closed set in $X^*$.
\end{prop}

$\Box$ \vspace{.5cm}

{\bfseries Example:} (F) Let $X$ be the disjoint union of $(0,1]$ and a countable set $A$.  Define a $T_1$, but not
Hausdorff, topology $\t$ on $X$
so that $U \in \t$ when $U \cap (0,1]$ is open and, in addition, either $U \cap A = \emptyset$ or else
$U \cap A$ is cofinite (i.e. has finite complement in $A$) and, in addition,
 $(0,a) \subset U \cap (0,1]$ for some $a > 0$.
If $Y = B \cup (0,1] \subset X$ with $B \subset A$
then $Y$ is compact if and only if $B \not= \emptyset$.  Choose $\{ B_n \}$ a strictly decreasing sequence of
cofinite subsets of $A$ with $B_0 = A$ and with empty intersection.
  We obtain a decreasing sequence $Y_n = B_n \cup (0,1]$
of compact sets with  intersection $Y = (0,1]$ noncompact. If $f(x) = 1 - 2^{-n}$ for $x \in B_n \setminus B_{n+1}$
and $f(x) = 1$ for $x \in (0,1]$, then $f$ is a fuzzy open set (not fuzzy closed), and for every $c \in (0,1), f_{[c}$ is compact.  By Lemma \ref{lemc02}, Condition L holds, but $f$ is not compact because $f_{[1}$ is not compact.

$\Box$ \vspace{1cm}

\bibliographystyle{amsplain}

\end{document}